\documentclass[11pt]{amsart}

\usepackage{amsmath,amsfonts,amsthm,amssymb,graphicx,tikz,color}
\usepackage[all,2cell,ps]{xy}

\theoremstyle{plain}
\newtheorem{thm}{Theorem}[section]

\newtheorem{cor}[thm]{Corollary}
\newtheorem{conjec}[thm]{Conjecture}
\newtheorem{qtn}[thm]{Question}
\newtheorem{prob}[thm]{Problem}

\theoremstyle{definition}
\newtheorem{defn}[thm]{Definition}

\newtheorem{rem}[thm]{Remark}

\theoremstyle{remark}

\DeclareMathOperator{\Aut}{Aut} 
 
\DeclareMathOperator{\SL}{SL} 
\DeclareMathOperator{\GL}{GL} 
 
 \DeclareMathOperator{\SO}{SO}

\DeclareMathOperator{\SU}{SU} 
\DeclareMathOperator{\SP}{Sp} \DeclareMathOperator{\Sp}{Sp}
 
 \DeclareMathOperator{\Comm}{Comm}

\DeclareMathOperator{\diag}{diag}

\DeclareMathOperator{\CAT}{CAT(0)}
\DeclareMathOperator{\Isom}{Isom}

\DeclareMathOperator{\Diff}{Diff}

\newcommand{\Ga}{\Gamma}

\newcommand{\Gam}{\Gamma}

\newcommand{\inv}{^{-1}}

\newcommand{\Fr}[1]{\ensuremath{\mathfrak{#1}}}

\newcommand{\fh}{\Fr h}

\newcommand{\Hy}{\ensuremath{\mathbb{H}}}
\newcommand{\CH}{\ensuremath{\mathbb{CH}}}
\newcommand{\HH}{\ensuremath{\mathbb{HH}}}
\newcommand{\CaH}{\ensuremath{\mathbb{OH}}}

\newcommand{\Q}{\ensuremath{\mathbb{Q}}}
\newcommand{\R}{\ensuremath{\mathbb{R}}}
\newcommand{\Z}{\ensuremath{\mathbb{Z}}}
\newcommand{\C}{\ensuremath{\mathbb{C}}}
\newcommand{\T}{\ensuremath{\mathbb{T}}}

\title[Superrigidity]{Superrigidity, arithmeticity, normal subgroups: results, ramifications and directions}
\author[D Fisher]{David Fisher}
\address{Department of Mathematics\\Indiana University\\Bloomington, IN 47405}
\email{fisherdm@indiana.edu}

\begin{document}

\begin{abstract}
This essay points to many of the interesting ramifications of Margulis' arithmeticity theorem, the superrigidity theorem, and normal subgroup theorem. We provide some history and background, but the main goal is to point to interesting open questions that stem directly or indirectly from Margulis' work and it's antecedents.
\end{abstract}

\maketitle



\section{Introduction}

 We begin with an informal overview of the events that inspire this essay and the work it describes.  For formal definitions and theorems, the reader will need to look into later sections of the paper, particularly Section \ref{section:margulistheorems}

In a few years in the early 1970's, Margulis transformed the study of lattices in semisimple Lie groups.  In this section and the next $G$ is a semisimple Lie group of real rank at least two with finite center and $\Gamma$ is an irreducible lattice in $G$. For brevity we will refer to these lattices as {\em higher rank lattices}. The reader new to the subject can always assume $G$ is $\SL(n,\R)$ with $n>2$.  We recall that a lattice is a discrete group where the volume of $G/\Gamma$ is finite, and that $\Gamma$ is called uniform if $G/\Gamma$ is compact and non-uniform otherwise.  In 1971, Margulis proved that non-uniform higher rank lattices are arithmetic, i.e. that they are commensurable to the integer points in some realization of $G$ as a matrix group \cite{MargulisArithmeticity}.  The proof used a result Margulis had proven slightly earlier on the non-divergence of unipotent orbits in the space $G/\Gamma$ \cite{Margulisnondivergence}.  This result on non-divergence of unipotent orbits has since played a fundamental role in homogeneous dynamics and its applications to number theory, a topic treated in many other essays in this volume. Margulis' arithmeticity theorem had been conjectured by Selberg and Piatetski-Shapiro. Piatetski-Shapiro had also conjectured the result on non-divergence of unipotent orbits \cite{Selberg}.  Both Selberg and Piateski-Shapiro had also conjectured the arithmeticity result for uniform lattices, but it was clear that that case requires a different proof, since the space $G/\Gamma$ is compact and questions of divergence of orbits do not make sense.

In 1974, Margulis resolved the arithmeticity question in truly surprising manner.  He proved his superrigidity theorem that classified the linear representations of a higher rank lattice $\Gamma$ over any local field of characteristic zero and used this understanding of linear representations to prove arithmeticity \cite{MargulisICM}.  Connections between arithmetic properties of lattices and the rigidity of their representations had been observed earlier by Selberg \cite{Selberg}. Important rigidity results had been proven in the local setting by Selberg, Weil, Calabi-Vesentini and others and in a more global setting by Mostow \cite{Selberg, Weil-I, Weil-II,Calabi, CalabiVesentini, Mostow-book}.  Despite this, the proof of the superrigidity
theorem and this avenue to proving arithmeticity were quite surprising at the time. The proof of the superrigidity theorem, though inspired by Mostow's study of boundary maps in his rigidity theorem, was also quite novel in the combination of ideas from ergodic theory and the study of algebraic groups.

Four years after proving his superrigidity and arithmeticity theorems, Margulis proved another remarkable theorem about higher rank lattices, the normal subgroup theorem. Margulis proofs of both superrigidity and the normal subgroup theorem were essentially dynamical and cemented ergodic theory as a central tool for studying discrete subgroups of Lie groups.

The main goal of this article is to give some narrative of the repercussions and echoes of Margulis' arithmeticity, superrigidity and normal subgroup theorems and the related results they have inspired in various areas of mathematics with some focus on open problems.  To keep true to the spirit of Margulis' work, some emphasis will be placed on connections to arithmeticity questions, but we will also feature some applications to settings where there is no well defined notion of arithmeticity.  For a history of the ideas that led up to the superrigidity theorem, we point the reader to survey written by Mostow at the time \cite{MostowSelberg} and to a discussion of history in another survey of the author \cite[Section 3]{FisherSurvey}.

In the next section of this essay we give precise statements of Margulis' results.  Afterwards we discuss various
more recent developments with an emphasis on open questions.  We do not attempt to give a totally comprehensive history.  In some cases, we mention results without giving full definitions and statements, simply in order to indicate the full breadth and impact of Margulis' results without ending up with an essay several times the length of the current one.  We mostly refrain from discussing proofs or only discuss them in outline. For a modern proof of superrigidity theorems, we refer the reader to the paper of Bader and Furman in this volume \cite{BaderFurman1}.
The proof is certainly along the lines of Margulis' original proof, but the presentation is particularly elegant
and streamlined.

\section{Arithmeticity and superrigidity: Margulis' results}
\label{section:margulistheorems}

For the purposes of this essay, we will always consider semisimple Lie groups with finite center
and use the fact that these groups can be realized as algebraic groups.  We will also have occasion
to mention algebraic groups over other local fields, but will keep the main focus on the case
of Lie groups for simplicity. Given a semisimple Lie group $G$, the real rank of $G$ is the dimension of the largest
subgroup of $G$ diagonalizable over $\R$.

Given an algebraic group $G$ defined over $\Q$, one can consider the integer points of the group, which
will denote by $G(\Z)$.  Arithmetic groups are a (slight) generalization of this construction.  We say two subgroups
$A_1$ and $A_2$ of $G$ are commensurable, if their intersection is finite index in each of them, i.e.  $[A_1 \cap A_2: A_i] < \infty$ for $i=1,2$.

A lattice $\Gamma <G$ is {\em arithmetic} if the following holds:  there is another semisimple algebraic Lie group $G'$
defined over $\Q$ with a homomorphism $\pi: G'\rightarrow G$ with $\ker(\pi)=K$ a compact group such that
$\Gamma$ is commensurable to $\pi(G(\Z))$.

A lattice $\Gamma$ in a product of groups $G_1 \times G_2$ is {\em irreducible} if the projection to each factor is indiscrete.
In most contexts this is equivalent to $\Gamma$ not being commensurable to a product of a lattice $\Gam_1$ in $G_1$ and a lattice $\Gam_2$ in $G_2$.
Irreducibility for a lattice in a product with more than two factors is defined similarly.  We can now state Margulis' arithmeticity theorem formally.

\begin{thm}[Margulis arithmeticity]
\label{thm:arithmeticity}
Let $G$ be a semisimple Lie group of real rank at least $2$ and $\Gamma <G$ an irreducible lattice, then $\Gamma$ is arithmetic.
\end{thm}

We will now state the superrigidity theorems and then briefly sketch the reduction of arithmeticity to superrigidity.  This requires considering representations over fields other than $\R$ or $\C$, namely representations finite extensions of the $p$-adic fields $\Q_p$. Together, these are all the local fields of characteristic zero. Superrigidity and arithmeticitiy are also known for groups over local fields of positive characteristic as both source and target by combined works of Margulis and Venkataramana and for targets groups over valued fields that are not necessarily local by the work of Bader and Furman in this volume \cite{MargulisBook, Venkataramana, BaderFurman1}.

To state the superrigidity theorem cleanly, we recall a definition.
Given a lattice $\Gamma<G$ and topological group $H$, we say a homomorphism $\rho:\Gamma \rightarrow H$ {\em almost extends to a homomorphism of $G$} if there are representations $\rho_G:G \rightarrow H$ and $\rho':\Gamma \rightarrow H$ such that $\rho_G$ is continuous, $\rho'(\Gamma)$ is pre-compact and commutes with $\rho_G(G)$, and such that $\rho(\gamma)=\rho_G(\gamma)\rho'(\gamma)$ for all $\gamma$ in $\Gamma$.  We can now state the strongest form of Margulis' superrigidity theorem that holds in our context:

\begin{thm}[Margulis Superrigidity]
\label{thm:superrigidity}
Let $G$ be a semisimple Lie group of real rank at least $2$, let $\Gamma <G$ be an irreducible lattice and $k$ a local field of characteristic zero. Then any homomorphism $\rho:\Gamma \rightarrow \GL(n,k)$ almost extends to a homomorphism of $G$.
\end{thm}

In many contexts this theorem is stated differently, with assumptions on the image of $\rho$. Assumptions often are chosen to allow $\rho$ to extend to $G$ rather than almost extend or to extend on a subgroup of finite index.  These assumptions are typically that $\rho(\Gamma)$ has simple Zariski closure and is not pre-compact, which guarantees extension on a finite index subgroup, and the additional assumption that the Zariski closure is center free, to guarantees an extension on all of $\Gamma$.  In many contexts where Margulis theorem is generalized beyond linear representations to homomorphisms to more general groups, it is only this type of special case which generalizes. The version we state here is essentially contained in \cite{MargulisBook}, at least when $G$ has finite center.  The case of infinite center is clarified in \cite{FisherMargulisLRC}. We will raise some related open questions later.

We sketch a proof of arithmeticity from superrigidity, for more details see e.g.  \cite[Chapter 6.1]{ZimmerBook} or \cite[Chapter IX]{MargulisBook}. First notice that since $\Gamma$ is finitely generated, the matrix entries of $\Gamma$ lie in a finitely generated field $k$ that is an extension of $\Q$.  Assume $G$ is simple and center free.  Then Theorem \ref{thm:superrigidity} implies every representation of $\Gamma$ either extends to $G$ or has bounded image. Note that $\Aut(\C)$ acts transitively on the set transcendentals numbers. So if we assume $k$ contains transcendentals, we can take the defining representation of $\Gamma$ and compose it with a sequence of automorphisms of $\C$ that send the trace of the image of some particular $\gamma$ to infinity.  It is obvious that this can't happen in a representation with bounded image and also not hard to check that it can't happen in one that extends to $G$.  This means that $k$ is a number field, so $\Gamma \subset G(k)$ and we want to show that $\Gamma \subset G(\mathcal{O}_k)$.  To see that $\Gamma \subset G(\mathcal{O}_k)$, assume not.  Then there is a prime $\mathfrak{p}$ of $k$ such that the image of $\Gamma$ in $G(k_{\mathfrak{p}})$ is unbounded where $k_{\mathfrak{p}}$ is the completion of $k$ for its $\mathfrak{p}$-adic valuation.  But this contradicts Theorem \ref{thm:superrigidity} since this unbounded representation should almost extend to $G$ with $\rho_G$ non-trivial and continuous and such $\rho_G$ cannot exist since $G(k_{\mathfrak{p}})$ is totally disconnected. To complete the proof, we want to show that $\Gamma$ is commensurable to $G(\mathcal{O}_k)$.  Assuming that $k$ is of minimal possible degree over $\Q$,we establish this by showing this amount to showing $G(\mathcal{O}_k)$ is already a lattice in $G$.  We show this by showing that for any Galois automorphism $\sigma$ of $k$ other than the identity, the map $\Gamma \rightarrow G(\sigma(\mathcal{O}_k))$ obtained by composing the identity with $\sigma$ has bounded image. This follows from the superrigidity theorem again simply because Galois conjugation does not extend to a continuous automorphism of the real points of $G$.

We mention next one additional application of the superrigidity theorem.  Let $V$ be a vector space, and assume a higher rank lattice $\Gamma$
acts on $V$ linearly.  A natural object of study with many applications is the cohomology of $\Gamma$ with coefficients in $V$.  The first cohomology is particularly useful for applications.  We have

\begin{thm}[Margulis first cohomology]
\label{theorem:margulish1}
  Let $\Gamma$ be a higher rank lattice and $V$ a vector space on which $\Gamma$ acts linearly, then $H^1(\Gamma, V)=0$.
\end{thm}

\noindent Let $H$ be the Zariski closure of $\Gamma$ in $\GL(V)$. The proof results from realizing that cocycles valued in $V$ correspond to representation into $H \ltimes V$, applying superrigidity to see that these representations all must be conjugate into $H$ and realizing that this implies the cocycle is trivial.  An important part of this argument is that we can apply superrigidity to the group $H \ltimes V$ which is neither semisimple nor reductive, since $V$ is contained in the unipotent radical.  We remark that if the image of $\Gamma$ in $\GL(V)$ is precompact and all simple factors of $G$ have higher rank, the result follows from Property $(T)$ for $\Gamma$. There are other ways of computing $H^1(\Gamma, V)$ using techniques from geometry and representation theory, but as far as the author knows, none of these quite recover the full statement of Theorem \ref{theorem:margulish1} in the case of non-uniform lattices see e.g. \cite{BorelWallach}.  These geometric and representation theoretic methods can also be used to show vanishing theorems concerning higher degree cohomology that are not accessible by Margulis' methods.

Margulis also proved a variant of superrigidity and arithmeticity for lattices with dense commensurators.  For a subgroup $\Gamma < G$ we define

$$\Comm_G(\Gamma)=\{ g \in G \text{ } | \text{ } g\Gamma g{\inv}\text{ and } \Gamma \text{ are commensurable}\}.$$

\noindent The next theorem was proved by Margulis at essentially the same time as the superrigidity theorem for higher rank lattices \cite{MargulisICM}.  The proof works independently of the rank of the ambient noncompact simple group $G$, but given Theorem \ref{thm:superrigidity} above, it is most interesting when the rank of $G$ is one.

\begin{thm}[Margulis Commensurator Superrigidity]
\label{thm:commsuperrigidity}
Let $G$ be a semisimple Lie group without compact factors, let $\Gamma <G$ be an irreducible lattice and let $\Lambda < \Comm_G(\Gamma)$ be dense in $G$ and $k$ a local field of characteristic zero. Then any homomorphism $\rho:\Lambda \rightarrow \GL(n,k)$ almost extends to a homomorphism of $G$.
\end{thm}

\noindent As before Margulis obtained a corollary concerning arithmeticity, that again is most interesting when the rank of $G$ is one.

\begin{cor}[Margulis Commensurator Arithmeticity]
\label{cor:arithmeticity}
Let $G$ be a semisimple Lie group, let $\Gamma <G$ be an irreducible lattice and let $\Lambda < \Comm_G(\Gamma)$ be dense in $G$, then $\Gamma$ is arithmetic.
\end{cor}

\noindent The argument that Theorem \ref{thm:commsuperrigidity} implies Corollary \ref{cor:arithmeticity} is essentially the same as the argument that Theorem \ref{thm:superrigidity} implies Theorem \ref{thm:arithmeticity}.
The converse to Corollary \ref{cor:arithmeticity}, that the communsurator of an arithmetic lattice is dense, was already known at the time of Margulis' work and is due to Borel \cite{BorelCrelle}.

An important related theorem of Margulis is the normal subgroup theorem.  We state here the version for
lattices in Lie groups \cite{MargulisNST}.

\begin{thm}
\label{thm:nst}
Let $G$ be a semisimple real Lie group of real rank at least $2$ and $\Gam<G$ an irreducible lattice.
Then any normal subgroup $N \lhd \Gam$ is either finite or finite index.
\end{thm}

\noindent One can view this statement as one about some kind of superrigidity of homomorphisms of
$\Gam$ to discrete groups:  either the representation is almost faithful or the image is bounded.
Knowing Theorem \ref{thm:arithmeticity}, Theorem \ref{thm:nst} can also be viewed as an arithmeticity theorem saying that any infinite normal subgroup of a higher rank arithmetic lattice is still an arithmetic lattice.
The proof of Theorem \ref{thm:nst} is quite different than the proof of Theorem \ref{thm:superrigidity}
but there is a longstanding desire to unify these phenomena in the context of higher rank lattices.









\section{Superrigidity and arithmeticity in rank one.}

The purpose of this section is to discuss lattices in rank $1$ simple Lie groups.  We discuss both known rigidity results and known constructions and raise some, mostly longstanding, questions.  The rank $1$ Lie groups are the isometry groups of various hyperbolic spaces:

\begin{enumerate}
  \item the group $\SO(n,1)$ is locally isomorphic to the isometry group of the $n$ dimensional hyperbolic space $\Hy^n$,
  \item the group $\SU(n,1)$ is locally isomorphic to the isometry group of the $n$ (complex) dimensional complex hyperbolic space $\CH^n$,
  \item the group $\SP(n,1)$ is locally isomorphic to the isometry group of the $n$ (quaternionic) dimensional quaternionic hyperbolic space $\HH^n$,
  \item the group $F_4^{-20}$ is the isometry group of the two dimensional Cayley hyperbolic plane $\CaH^2$.
\end{enumerate}

\noindent Exceptional isogenies between Lie groups yield isometries between some low dimensional hyperbolic spaces, namely that $\Hy^2=\CH^1$, that $\HH^1=\Hy^4$ and that $\CaH^1=\Hy^8$.

The strongest superrigidity and arithmeticity results for rank one groups generalize Margulis' results completely to lattices in $\Sp(n,1)$ and  $F_4^{-20}$.  There are also numerous interesting partial results for lattices in the other two families of rank one Lie groups $\SO(n,1)$ and $\SU(n,1)$.

At the time of Margulis' proof of arithmeticity, non-arithmetic lattices were only known to exist in $\SO(n,1)$ when $2 \leq k \leq 5$.  No non-arithmetic lattices were known in the other rank one simple groups.  Margulis asked about the other cases in \cite{MargulisICM}. In this section we will also discuss known results, including other criteria for arithmeticity of lattices in rank $1$ groups and known examples of non-arithmetic lattices.

\subsection{Quaternionic and Cayley Hyperbolic spaces}
\label{sec:sp}

In this subsection we describe the developments that proved that all lattices in $\SP(n,1)$ for $n>1$ and $F_4^{-20}$ are arithmetic. The first major result in this direction, concerning rigidity of quaternionic and Cayley hyperbolic lattices was proved by Corlette \cite{Corlette-Annals}.  He showed

\begin{thm}[Corlette]
\label{thm:corlette}
Let $G= \Sp(n,1)$ for $n>1$ or $G=F_4^{-20}$ and $\Gamma<G$ be a lattice.  Let $H$ be a real simple Lie group with finite center
and $\rho:\Gamma \rightarrow H$ a homomorphism with unbounded Zariski dense image.  Then $\rho$ almost extends to $G$.
\end{thm}

\noindent {\bf Remarks:}
\begin{enumerate}
  \item When $n=1$, the group $\Sp(1,1)$ is isomorphic to $\SO(4,1)$.
  \item In this setting, one can replace that $\rho$ almost extends with the statement that $\rho$ extends on a subgroup  of finite index.
\end{enumerate}

\noindent The proof of Corlette's theorem has two main steps.  The first is the existence of a $\Gamma$ equivariant harmonic map from $\HH^n$ or $\CaH^2$ to $H/K$, the symmetric space associated to $H$. This step is contained in
earlier work of Corlette or Donaldson, see also Labourie \cite{CorletteJDG, Donaldson, Labourie}.  Corlette then proves a Bochner formula that allows him to conclude the harmonic map is totally geodesic, from which the result follows relatively easily.  This work is inpsired by earlier work of Siu that proved generalizations of Mostow rigidity using harmonic map techniques \cite{SiuMostow}.  The idea of using harmonic maps to prove superrigidity theorems was well known at the time of Corlette's work and is often attributed to Calabi.

Following Corlette's work, Gromov and Schoen developed the existence and regularity theory of harmonic maps to buildings in order to prove \cite{GromovSchoen}:

\begin{thm}[Gromov-Schoen]
Let $G= \Sp(n,1)$ for $n>1$ or $G=F_4^{-20}$ and $\Gamma<G$ be a lattice.  Let $H$ be a simple algebraic group over a non-Archimedean local field with finite center and $\rho:\Gamma \rightarrow H$ a homomorphism with Zariski dense image.  Then $\rho$ has bounded image.
\end{thm}

\noindent The main novelty in the work of Gromov and Schoen is to prove existence of a harmonic map into certain singular spaces with enough regularity of the harmonic map to apply Corlette's Bochner inequality argument.  The harmonic map is to the Euclidean building associated to $H$ by Bruhat and Tits \cite{BruhatTits}, and it is easy to see that there are no totally geodesic maps from hyperbolic spaces to Euclidean buildings.

Combining these two results with arguments of Margulis' deduction of arithmeticity from superrigidity, we can deduce:

\begin{thm}
\label{thm:sparithmeticity}
Let $G= \Sp(n,1)$ for $n>1$ or $G=F_4^{-20}$ and $\Gamma<G$ be a lattice, then $\Gamma$ is arithmetic.
\end{thm}

We mention here a related result of Bass-Lubotzky that answered a question of Platonov \cite{BassLubotzky, LubotzkySuperrigid2}.  Namely Platonov asked if any linear group that satisfied the conclusion of the superrigidity theorem was necessarily an arithmetic lattice.  Bass and Lubotzky produce counter-examples as subgroups $\Delta < \Gamma \times \Gamma$ such that $\diag(\Gamma) < \Delta$ where $\Gamma <G$ is a lattice and $G$ is either $F_4^{-20}$ or $\SP(n,1)$ for $n>1$. The proofs involve a number of new ideas but depend pivotally on the work of Corlette and Gromov-Schoen to prove the required superrigidity results. In the examples produced by Bass and Lubotzky, the proof that $\Delta$ is superrigid is always deduced from the known superrigidity of $\diag(\Gamma)$.  The fact that $\Gamma$ is  a hyperbolic group in the sense of Gromov plays a key role in constructing $\Delta$.

\begin{qtn}
\label{qtnsuperrigid}
Are there other superrigid non-lattices?  Can one find a superrigid non-lattice that is Zariski dense in higher rank simple Lie group?  Can one find a superrigid non-lattice that does not contain a superrigid lattice?  A superrigid non-lattice which is a discrete subgroup of a simple non-compact Lie group?
\end{qtn}

\subsection{Results in real and complex hyperbolic geometry}

\subsubsection{Non-arithmetic lattices: constructions and questions}
\label{subsub:constructionsandquestions}
To begin this subsection I will discuss the known construction of non-arithmetic lattices in $\SO(n,1)$ and $\SU(n,1)$. To begin slightly out of order, I emphasize one of the most important open problems in the area, borrowing wording from Margulis in \cite{MargulisPCR}.

\begin{qtn}
\label{qtn:nonarithmeticSU}
For what values of $n$ does there exist a non-arithmetic lattice in $\SU(n,1)$?
\end{qtn}

The answer is known to include $2$ and $3$.  The first examples were constructed by Mostow in \cite{MostowPJM} using reflection group techniques. The list was slightly expanded by Mostow and Deligne using monodromy of hypergeometric functions \cite{DeligneMostow, MostowGPL}.  The exact same list of examples was rediscovered/re-interpreted by Thurston in terms of conical flat structures on the $2$ sphere \cite{ThurstonShapes}, see also \cite{SchwartzNotesonShapes}. There is an additional approach via algebraic geometry suggested by Hirzebruch and developed by him in collaboration with Barthels and H\"{o}fer \cite{BHH}. More examples have been discovered recently by Couwenberg, Heckman, and Looijenga using the Hirzebruch style techniques and by Deraux, Parker and Paupert using complex reflection group techniques \cite{CHL, DPP, DPP1, Deraux3d}.  But as of this writing there are only $22$ commensurability classes of non-arithmetic lattices known in $\SU(2,1)$ and only $2$ known in $\SU(3,1)$.   An obvious refinement of Question \ref{qtn:nonarithmeticSU} is

\begin{qtn}
\label{qtn:nonarithmeticSUcomm}
For what values of $n$ do there exist infinitely many commensurablity classes of non-arithmetic lattice in $\SU(n,1)$?
\end{qtn}

\noindent We remark here that the approach via conical flat structures was extended by Veech and studied further by Ghazouani and Pirio \cite{VeechFlatSurf, GhazouaniPirio1}.  Regretably this approach does not yield more non-arithmetic examples.  It seems that the reach of this approach is roughly equivalent to the reach of the approach via monodromy of hypergeometric functions, see \cite{GP2}.  There appears to be some consensus among experts is that the answer to both Question \ref{qtn:nonarithmeticSU} and Question \ref{qtn:nonarithmeticSUcomm} should be ``for all $n$", see e.g. \cite[Conjecture 10.8]{Kapovich}.  We note here that Margulis' own wording as used above is more guarded.

At the time of Margulis work the only known non-arithmetic lattices in $\SO(n,1)$ for $n>2$ were constructed by Makarov and Vinberg by reflection group methods \cite{Makarov, VinbergFirst}.  It is known  by work of Vinberg that these methods will only produce non-arithmetic lattices in dimension less than $30$ \cite{Vinberg}.  The largest known non-arithmetic lattice produced by these methods is in dimension $18$ by Vinberg and the full limits of reflection group constructions is not well understood \cite{VinbergHigh}. We refer the reader to \cite{Belolipetsky} for a detailed survey. The following question seems natural:

\begin{qtn}
In what dimensions do there exist lattices in $\SO(n,1)$ or $\SU(n,1)$ that are commensurable to non-arithmetic reflection groups?  In what dimensions do there exist lattices in $\SO(n,1)$ or $\SU(n,1)$ that are commensurable to arithmetic reflection groups?
\end{qtn}

For the real hyperbolic setting, there are known upper bounds of $30$ for non-arithmetic lattices and $997$ for any lattices. The upper bound of $30$ also applies for arithmetic uniform hyperbolic lattices \cite{Vinberg, Belolipetsky}.  In the complex hyperbolic setting, there seem to be no known upper bounds, but a similar question recently appeared in e.g. \cite[Question 10.10]{Kapovich}.  For a much more detailed survey of reflection groups in hyperbolic spaces, see \cite{Belolipetsky}.

A dramatic result of Gromov and Piatetski-Shapiro vastly increased our stock of non-arithmetic lattices in $\SO(n,1)$
by an entirely new technique:

\begin{thm}[Gromov and Piatetski-Shapiro]
For each $n$ there exist infinitely many commensurability classes of non-arithmetic uniform and non-uniform lattices in $\SO(n,1)$.
\end{thm}

The construction in \cite{GPS} involves building hybrids of two arithmetic manifolds by cutting and pasting along totally geodesic codimension one submanifolds.  The key observation is that non-commensurable arithmetic manifolds can contain isometric totally geodesic codimension one submanifolds. This method has been extended and explored by many authors for a variety of purposes, see for example \cite{Agol, BelThom, SevenSamuraiShort, GelanderLevit}.  It has also been proposed that one might build non-arithmetic complex hyperbolic lattices using a variant of this method, though that proposal has largely been stymied by the lack of codimension one totally geodesic codimension one submanifolds.
The absence of codimension $1$ submanifolds makes it difficult to show that attempted ``hybrid" constructions yield discrete groups.  For more information see e.g. \cite{Paupert, PaupertWells, Wells} and
\cite[Conjecture 10.9]{Kapovich}.  We point out here that the results of Esnault and Groechenig discussed below as
Theorem \ref{EGcorollary} implies that the ``inbreeding" variant of Agol and Belolipetsky-Thomsen \cite{Agol, BelThom} cannot produce non-arithmetic manifolds in the complex hyperbolic setting even if the original method of Gromov and Piatetski-Shapiro does.

In \cite{GPS}, Gromov and Piatetski-Shapiro ask the following intriguing question:

\begin{qtn}
\label{qtn:pieces}
Is it true that, in high enough dimensions, all lattices in $\SO(n,1)$ are built from sub-arithmetic pieces?
\end{qtn}

\noindent
The question is somewhat vague, and sub-arithmetic is not defined in \cite{GPS}, but a more precise starting point is:

\begin{qtn}
\label{question:tg}
For $n>3$, is it true that any non-arithmetic lattice in $\Gamma< \SO(n,1)$ intersects some conjugate of $\SO(n-1,1)$ in a lattice?
\end{qtn}

\noindent This is equivalent asking if every finite volume non-arithmetic hyperbolic manifold contains a closed codimension one totally geodesic submanifold.   Both reflection group constructions and hybrid constructions all contain such submanifolds.   It seems the consensus in the field is that the answer to this question should be no, but we know of no solid evidence for that belief.  It is also not known to what extent the hybrid constructions and reflection group constructions build distinct examples. Some first results, indicating that the classes are different, are contained in \cite[Theorem 1.7]{FLMS} and in \cite[Theorem 1.5]{Mila}.

It is worth mentioning that our understanding of lattices in $\SO(2,1)$ and $\SO(3,1)$ is both more developed and very different.
Lattices in $\SO(2,1)$ are completely classified, but there are many of them, with the typical isomorphism class of lattices
having many non-conjugate realizations as lattices, parameterized by moduli space.  In $\SO(3,1)$, Mostow rigidity means there
are no moduli spaces.  But Thurston-Jorgensen hyperbolic Dehn surgery still allows one to construct many ``more" examples of lattices,
including ones that yield a negative answer to Question \ref{question:tg}.  There remains an interesting sense in which the answer
to Question \ref{qtn:pieces} could still be yes even for dimension $3$

\begin{qtn}
\label{qtn:dehnsurgery}
Can every finite volume hyperbolic $3$-manifold be obtained as Dehn surgery on an arithmetic manifold?
\end{qtn}

\noindent To clarify the question, it is known that every finite volume hyperbolic $3$-manifold is obtained as a topological manifold  by Dehn surgery on some cover of the figure $8$ knot complement, which is known to be the only arithmetic knot complement \cite{HLM, ReidKnot}.  What is not known is whether one can obtain the geometric structure on the resulting three manifold as geometric deformation of the complete geometric structure on the arithmetic manifold on which one performs Dehn surgery.

\subsubsection{Arithmeticiy, superrigidity and totally geodesic submanifolds}

This section concerns recent results by Bader, the author, Miller and Stover, motivated by questions of McMullen and Reid in the case of real hyperbolic manifolds.  Throughout this section a geodesic submanifold will mean a closed immersed, totally geodesic submanifold. (In fact all results can be stated also for orbifolds but we ignore this technicality here.) A geodesic submanifold {\emph{maximal} if it is not contained in a proper geodesic submanifold of smaller codimension.

For arithmetic manifolds, the presence of one maximal geodesic submanifold can be seen to imply the existence of infinitely many.  The argument involves lifting the submanifold $S$ to a a finite  cover $\tilde M$ where an element $\lambda$ of the commensurator acts as an isometry.  It is easy to check that $\lambda(S)$ can be pushed back down to a geodesic submanifold of $M$ that is distinct from $S$.  This was perhaps first made precise in dimension $3$ by Maclachlan--Reid and Reid \cite{MRTG, ReidTG}, who also exhibited the first hyperbolic $3$-manifolds with no totally geodesic surfaces.

In the real hyperbolic setting the main result from \cite{BFMS} is

\begin{thm}[Bader, Fisher, Miller, Stover]
\label{thm:bfms1}
Let $\Gam$ be a lattice in $\SO_0(n, 1)$. If the associated locally symmetric space contains infinitely many maximal geodesic submanifold of dimension at least $2$, then $\Gam$ is arithmetic.
\end{thm}

\begin{rem} \begin{enumerate}

\item The proof of this result involves proving a superrigidity theorem for \emph{certain} representations of the lattice in $\SO(n,1)$.  As the conditions required become a bit technical, we refer the interested reader to
    \cite{BFMS}.  The superrigidity is proven in the language introduced in \cite{BaderFurman1}.

\item At about the same time, Margulis and Mohammadi gave a different proof for the case $n=3$ and $\Gam$ cocompact \cite{MM}.  They also proved  a superrigidity theorem, but both the statement and the proof are quite different than in \cite{BFMS}.

\item A special case of this result was obtained a year earlier by the author, Lafont, Miller and Stover \cite{FLMS}. There we prove that  a large class of non-arithmetic manifolds have only finitely many maximal totally geodesic submanifolds.  This includes all the manifolds constructed by Gromov and Piatetski-Shapiro but not the examples constructed by Agol and Belilopetsky-Thomsen.

\end{enumerate}
\end{rem}

In the context of Margulis' work it is certainly worth mentioning that Theorem \ref{thm:bfms1} has a reformulation entirely in terms of homogeneous dynamics and that homogenenous dynamics play a key role in the proof.  It is also interesting that a key role is also played by dynamics that are not quite homogeneous but that take place on a projective bundle over the homogeneous space $G/\Gamma$.

Even more recently the same authors have extended this result to cover the case of complex hyperbolic manifolds.

\begin{thm}[Bader, Fisher, Miller, Stover] \label{thm:charithmetic}
Let $n \ge 2$ and $\Gam < \SU(n,1)$ be a lattice and $M=\CH^n/\Gamma$. Suppose that $M$ contains infinitely many maximal totally geodesic submanifolds of dimension at least $2$. Then $\Gam$ is arithmetic.
\end{thm}

As before this is proven using homogeneous dynamics, dynamics on a projective bundle over $G/\Gamma$, and a superrigidity theorem.  Here the superrigidity theorem is even more complicated than before and depends also on results of Simpson and Pozzetti \cite{Simpson, Pozzetti}.

The results in this section provide new evidence that totally geodesic manifolds play a very special role in non-arithmetic lattices and perhaps provide some evidence that the conventional wisdom on Questions \ref{question:tg} and \ref{qtn:nonarithmeticSU} should be reconsidered.

\subsubsection{Other superrigidity and arithmeticity results for lattices in $\SO(n,1)$ and $\SU(n,1)$}

The combination of the results in the last section and Margulis' commensurator superrigidity theorem, as well
as questions in \ref{subsub:constructionsandquestions} raise the following:

\begin{qtn}
\label{qtn:rankonesuperrigid}
Let $\Gamma <G$ be a lattice where $G=\SO(n,1)$ or $\SU(n,1)$.  What conditions on a representation $\rho: \Gamma \rightarrow \GL(m,k)$ implies that $\rho$ extends or almost extends?  What conditions on $\Gam$ imply that $\Gam$ is arithmetic?
\end{qtn}

\noindent For $\SU(n,1)$ Margulis asks a similar, but more restricted, question in \cite{MargulisPCR}.  He asks whether there might be particular lattices in $\SU(n,1)$ where superrigidity holds without restrictions on $\rho$ as in the higher rank, quaternionic hyperbolic and Cayley hyperbolic cases.

A very first remark is that for many $\Gam$ as above it is known that there are surjections of $\Gam$ on both abelian and non-abelian free groups. This suggests that one might want to study faithful representations or ones with finite kernel, though surprisingly very few known superrigidity results explicitly assume faithfulness of the representation.
The main counter-example to this is the following theorem of Shalom \cite{Shalom-Annals}.  We recall that for a discrete group $\Delta$ of a rank one simple Lie group, $\delta(\Delta)$ is the Hausdorff dimension of the limit set of $\Delta$.  The limit set admits many equivalent definitions see e.g. \cite{Shalom-Annals} for discussion.

\begin{thm}[Shalom]
\label{thm:shalom}
Let $\Gam<G$ be a lattice where $G=\SO(n,1)$ or $\SU(n,1)$.  Let $\rho: \Gam \rightarrow H$ be a discrete, faithful representation where $H$ is either $\SO(m,1)$ or $\SU(m,1)$.  Then $\delta(\Gam) \leq \delta(\rho(\Gam))$.
\end{thm}

\noindent  Shalom actually proves a result for non-faithful discrete representations as well, relating the dimension of the limit set of the image and the kernel to the dimension of the limit set of the lattice.  Shortly after Shalom proved the above theorem, Besson, Courtois and Gallot proved that equality only occurs in the case where the representation almost extends \cite{BCGSchwartzLemma}.  The methods of Besson, Courtois and Gallot, the so-called barycenter mapping, have been used in many contexts.  The key ingredient in Shalom's proofs, understanding precise decay rates of matrix coefficients, has not been exploited nearly as thoroughly for applications to rigidity.  For either Shalom's techniques or the barycenter map technique, the utility of the methods are currently limited by the requirement that the representation have discrete image.

Relatively few other superrigidity or arithmeticity type results are known for real hyperbolic manifolds but a plethora of other interesting phenomena have been discovered in the complex hyperbolic setting.  We begin with some of the most recent, which involves a bit of detour in a surprising direction.

Simpson's work on Higgs bundles and local systems focuses broadly on the representation theory of $\pi_1(M)$ where $M$ is a complex projetive variety or more generally a complex quasi-projective variety \cite{Simpson}.  This is related to our concerns because when $G=\SU(n,1)$ then $M= K \backslash G /\Gam$ is a projective variety when $\Gam$ is compact and quasi-projective when it is not. We say a representation $\rho: \Gam \rightarrow H$ is \emph{rigid} or \emph{infinitesimally rigid} if the first cohomology $H^1(\Gam, \fh)$ vanishes where $\fh$ is the Lie algebra of $H$. For $G=\SU(n,1), H=G$ and $\rho$ the defining representation $\rho :\Gam \rightarrow H$, vanishing of this cohomology group is a result of Calabi-Vesentini \cite{CalabiVesentini}. We state Simpson's main conjecture only in the projective case to avoid technicalities \cite{Simpson}:

\begin{conjec}
\label{conjecture:simpson}
Let $M$ be a projective variety and $\rho: \pi_1(M) \rightarrow \SL(n, \C)$ an infinitesimally rigid representation.  Then $\rho(\Gamma)$ is integral, i.e. there is a number field $k$ such that $\rho(\pi_1(M))$ is contained in the integer points $\SL(n, \mathcal{O}_k)$.
\end{conjec}

We state the conjecture for $\SL$ targets rather than $\GL$ targets to avoid a technical finite determinant
condition. We remark that higher rank irreducible K\"{a}hler locally symmetric spaces of finite volume provide examples where Simpson's conjecture follows from Margulis' arithmeticity theorem. Recent work of Esnault and Groechenig prove this result in many cases \cite{EG1, EsnaultGroechenig}. In particular their results have the following as a (very) special case:

\begin{thm}[Esnault-Groechenig]
\label{EGcorollary}
Let $\Gam < \SU(n,1)$ be a lattice with $n>1$, then $\Gam$ is integral.  I.e. there is a number field $k$ and $k$ structure on $\SU(n,1)$ such that $\Gam < \SU(n,1)(\mathcal{O}_k)$.
\end{thm}

The theorem is immediate from the results in \cite{EsnaultGroechenig} for the case of cocompact lattices.  For an explanation of how it also follows in the noncocompact case see \cite{BFMS2}.  We also note that a construction of Agol as extended by Belilopetsky-Thomson shows that the analogous result fails in $\SO(n,1)$ \cite{Agol, BelThom}.  I.e. there are non-integral lattices, both cocompact and noncocompact, in $\SO(n,1)$ for every $n$.

We note here that the proof Esnault and Groechenig does not pass through a superrigidity theorem.  In the context of this paper, one might expect this, but the methods of \cite{EsnaultGroechenig} depend on algebraic geometry and deep results of Lafforgue on the Langlands program \cite{Lafforgue}.  However in this context one might also ask the following:

\begin{qtn}
\label{qtn:sun1nonarch}
Let $\Gam < \SU(n,1)$ be a lattice and $n>1$.  Assume $k$ is a totally disconnected local field, $H$ is a simple algebraic group over $k$ and  $\rho: \Gam \rightarrow H$ is a Zariski dense, faithful representation.  Is $\rho(\Gam)$
compact?
\end{qtn}

We mention here a question from our paper with Larsen, Stover and Spatzier \cite{FLSS} that aims at understanding the degree to which a lattice in $\SO(n,1)$ can fail to be integral by studying the $p$-adic representation theory of these groups.

\begin{qtn}\label{qtn:DFRepn}
Let $\Sigma_g$ be a surface group of genus $g \ge 2$. Is there a discrete and faithful representation of $\Sigma_g$ into $\Aut(Y)$ for $Y$ a locally compact Euclidean building? Can we take $Y$ to be a finite product of bounded valence trees?
\end{qtn}

More generally one can ask the same questions with $\Sigma_g$ replaced by a lattice $\Gam$ in $G=\SO(n,1)$.  Once $n>2$, it is known that $\Gam$ is contained in the $k$ points of $G$ for some number field $k$.  To understand the extent to which $\Gam$ fails to be integral, it suffices to consider the case where $Y$ is the building associated to some $p$-adic group $G(k_p)$.

There is one other context in which enough superrigidity results are known to imply arithmeticity, namely Klingler's work on fake projective planes \cite{KlinglerFPP}.

\begin{defn}
A fake projective plane is a complex projective surface with the same Betti numbers as $P(\C^2)$ that is not
biholomorphic to $P(\C^2)$.
\end{defn}

 Results of Yau on the Calabi conjecture show that any fake projective plane is of the form $\CH^2/\Gamma$ with $\Gamma$ a cocompact lattice \cite{YauCalabi}. Let $G=\SU(2,1)$, we can further assume that $K \backslash G/ \Gam =M$ satisfies the condition that $c_1^2 = 3 c_2 = 9$ where $c_1$ and $c_2$ are the first and second Chern numbers of $M$. Yau's work implies complex ball quotients satisfying these conditions are exactly the fake projective planes.  Klingler then shows that $\Gam$ is arithmetic

\begin{thm}[Klingler]
\label{klinglerfpp}
If $M$ is a fake projective plane, then $\Gam = \pi_1(M)$ is arithmetic.
\end{thm}

\noindent  This result is striking since the condition for arithmeticity is purely topological. The proof uses superrigidity theorems proven using harmonic map techniques as in Subsection \ref{sec:sp}. Following Klingler's work, the fake projective planes were classified and further studied by Prasad-Yeung and Cartwright-Steger \cite{PrasadYeung, CartwrightSteger}.  There turn out to be exactly $50$ examples.  We note that this is precisely fifty and not fifty up to commensurability and that some of these examples are commensurable.  The topological condition of being a fake projective plane is not invariant under passage to finite covers.

Two more recent result of Klingler and collaborators are also intriguing in this context. In the first of these papers he shows that for certain lattices $\Gam$ in $\SU(n,1)$, the representation theory of $\Gam$ is very restricted as long as one considers representations in dimension below $n-1$ \cite{KlinglerSym}.  The results there are proven by showing that holomorphic  symmetric differentials control the linear representation theory of fundamental groups of compact K\"{a}hler manifolds.  In a later paper by Brunebarbe, Klingler and Totaro, the authors extend this to investigate the case of compact K\"{a}hler manifolds without holomorphic symmetric differentials  \cite{BKT}.

A different direction for the study of representations of complex hyperbolic lattices was introduced by Burger and Iozzi in \cite{BurgerIozzi}.  They introduce a notion of a maximal representation of a lattice $\Gam$ in $G=\SU(n,1)$
generalizing a definition of Toledo in the case of $\SU(1,1) \cong \SL(2,\R)$ \cite{Toledo}.  Burger-Iozzi show that
maximal representation of $\Gamma$ into $\SU(m,1)$ extend to $G$.  The proof uses a result on incidence geometry generalizing an earlier result of Cartan to the measurable setting \cite{Cartan}.  The definitions and results of this paper were further extended to the case of $\SU(p,q)$ targets when $p \neq q$ by Pozzetti in her thesis \cite{Pozzetti}.  The proof of Theorem \ref{thm:charithmetic} uses Pozzetti's version of the Cartan theorem.  More recently Koziarz and Maubon extended the result to include the case where $p=q$ and reproved all earlier results using techniques of harmonic maps and Higgs bundles \cite{KoziarzMaubon}.

\section{Orbit equivalence rigidity}
\label{section:oe}

This section mostly serves to point to a broad area of research that I will not attempt to summarize or survey in any depth.

\begin{defn}
\label{def:oe}
Let $(S,\mu)$ be a finite measure space with an ergodic $G$ action and $(S',\mu)$ a finite measure space with an ergodic $G'$ action.  We say the actions are \emph{orbit equivalent} if there are conull Borel sets $S_0 \subset S$
and $S_0' \subset S'$ and a measure class preserving isomorphism $\phi: S_0 \rightarrow S_0'$ such that $s$ and $t$ are in the same $G$ orbit if and only if $\phi(s)$ and $\phi(t)$ are in the same $G'$ orbit.
\end{defn}

In a remarkable result in \cite{ZimmerOE}, Zimmer further developed the ideas in Margulis' proof of superrigidity
to prove:

\begin{thm}
\label{Zimmer:oe}
Let $G_1$ and $G_2$ be center-free connected simple Lie groups and assume $\R$-rank$(G_1)>1$.  Let $(S_i,\mu_i)$ be probability measure spaces with ergodic actions of $G_i$ for $i=1,2$.  If the actions are orbit equivalent, then they are conjugate.
\end{thm}

The key ingredient in the proof of Theorem \ref{Zimmer:oe} is Zimmer's cocycle superrigidity theorem.  We do not state this here but point the reader to \cite{ZimmerBook, FisherMargulisLRC} for detailed discussions.

An important further development in the theory comes in work of Furman, who extends Zimmer's results on orbit equivalence to lattices \cite{FurmanGromov, Furman}.  We do not give a comprehensive discussion but state one result.

\begin{thm}
\label{furman:oe}
Let $G$ be a center-free connected simple Lie group and assume $\R$-rank$(G)>1$. Let $\Gam_1<G$ be a lattice and let $\Gam_2$ be any finitely generated group.  Let $(S_i,\mu_i)$ be probability measure spaces with ergodic actions of $\Gam_i$ for $i=1,2$.  If the actions of $\Gamma_i$ on $(S,\mu_i)$  are orbit equivalent, then $\Gam_2$ is virtually a lattice in $G$.
\end{thm}

\noindent Here virtually means there is a finite index subgroup of $\Gam_2$ whose quotient by a finite normal subgroup is a lattice in $G$. Furman also shows that there is a unique obstruction to conjugacy of the actions.

Following these results, the study of orbit equivalence rigidity became a rich topic in which many rigidity results are known, many of which depend on cocycle superrigidity theorems.  We do not attempt a survey but point to one written earlier by Furman \cite{FurmanSurvey}.

\section{The Zimmer program}
\label{section:Zimmerprogram}

In 1983, Zimmer proposed a  number of conjectures about actions of higher rank simple Lie groups and their lattices on compact manifolds \cite{ZimmerBulletin,ZimmerICM}.  These conjectures were motivated by a number of Zimmer's own theorems, including the cocycle superrigidity theorem mentioned in the last section.  But perhaps the clearest motivation is as a \emph{non-linear} analogue of Margulis' superrigidity theorem.  These conjectures have led to a tremendous amount of activity, see the author's earlier survey and recent update \cite{FisherSurvey, FisherSurvey2} for more information.  Here we focus only on two aspects: the recent breakthrough made by Brown, the author and Hurtado, and a statement of a general conjectural superrigidity theorem for $\Diff(M)$ targets.

The clearest conjecture made by Zimmer predicted that any action of a higher rank lattice on a compact manifold of sufficiently small dimension should preserve a Riemannian metric.  Since the isometry group of a compact manifold is a compact Lie group, this, together with Margulis' superrigidity theorem, often implies the action factors through a finite quotient of the lattice.  The recent work of Brown, the author and Hurtado makes dramatic progress on this conjecture and completely resolves it in several key cases \cite{BFH, BFH2, BFH3}.  For example we have:

\begin{thm}
\label{thm:bfh}[Brown, Fisher, Hurtado]
Let $\Gamma$ be a lattice in $\SL(n,\R)$, let $M$ be
a compact manifold and let $\rho:\Ga \rightarrow \Diff(M)$ be a homomorphism.
Then
\begin{enumerate}
  \item if $\dim(M) < n-1$, the image of $\rho$ is finite;
  \item if $\dim(M) < n$ and $\rho(\Gamma)$ preserves a volume form on $M$,
  then the image of $\rho$ is finite.
\end{enumerate}
\end{thm}

\noindent This result is sharp, since $\SL(n,\R)$ acts on the projective space $P(\R^n)$ and $\SL(n,\Z)$ acts on the torus $\T^n$.
The papers with Brown and Hurtado prove results about all lattices in all simple Lie groups $G$ of higher real rank;
but are only sharp for certain choices of $G$.  In particular results about volume preserving actions are only sharp
for $\SL(n,\R)$ and $\Sp(2n, \R)$ while results about actions not assumed to preserve volume are sharp for all split simple groups.  See \cite{Cantat, BFH3} for more discussion.

The most naive version of the Zimmer program is perhaps the following

\begin{qtn}
Let $G$ be a simple Lie group of higher real rank, $\Gamma <G$ a lattice and $M$ a compact manifold.
Can one understand all homomorphisms $\rho : \Gamma \rightarrow \Diff(M)$?  If $\omega$ is a volume
form on $M$, can one classify all homomorphisms $\rho : \Gamma \rightarrow \Diff(M, \omega)$?
\end{qtn}

\noindent The careful reader will notice a slight variation in wording in the two questions.  This is due
to the fact that non-volume preserving actions are known to be non-classifiable.  In particular the parabolic
induction described by Stuck in \cite{Stuck} shows that even homomorphisms $\rho: G \rightarrow \Diff(M)$
cannot be classified.  In particular Stuck shows that given two vector fields $X$ and $Y$ on a compact manifold $M$
and a parabolic subgroup $Q$ in $G$, one can construct two homomorphisms $\rho_X, \rho_Y: G \rightarrow \Diff((G \times M)/Q)$ such that $\rho_X$ and $\rho_Y$ are conjugate if and only if the flows generated by $X$ and $Y$ on $M$ are conjugate.

We briefly describe Stuck's construction.  Any parabolic subgroup $Q <G$ admits a homomorphism $\phi: Q \rightarrow \R$. Any vector field $X$ on $M$ defines an $\R$ action which we denote by $\bar{\rho}_X: \R \times M \rightarrow M$.
We define a $Q$ action on $G \times M$ by $(g,m)q= (gq\inv, \bar{\rho}(\phi(q)))$.  As this commutes with the left
$G$ action on the first variable, we obtain an action $\rho_X$ of $G$ on $(G\times M)/Q$.  The space $(G \times M)/Q$ is a manifold and in fact an $M$ fiber bundle over $G/Q$.  It is transparent that applying the construction to two vector fields $X$ and $Y$ on manifolds $M$ and $M'$, the $G$ actions are conjugate if and only $\bar{\rho}_X$ and $\bar{\rho}_Y$ are.  The following seems accessible

\begin{prob}
\label{problem:stuck}
If $\rho_X$ and $\rho_Y$ are conjugate as $\Gam$ actions, are $\bar{\rho}_X$ and $\bar{\rho}_Y$ conjugate as $\R$ actions?
\end{prob}

The main goal of this section is to describe a conjectural picture of all $\Gam$-actions on compact manifolds $M$ in terms of $G$ actions.  This is very much in the spirit of Margulis' superrigidity theorem and
we begin by restating that theorem slightly.  Given a higher rank lattice $\Gam$ in Lie group $G$ there is a natural
compact extension $G \times K$ of $G$ in which $\Gamma$ sits diagonally as a lattice.    Here $K$ is a product of two groups $K=K_1 \times K_2$ where $K_1$ is totally disconnected and $K_2$ is a Lie group.  The group $K_1$ is the profinite completion of $\Gam$.  The group $K_2$ is the compact Lie group such that $\Gam$ is commensurable to the integral points in $G' = G \times K_2$. The definition of arithmeticiy in Section \ref{section:margulistheorems} ensures that $K_2$ exists. We note here that $K_2$ is semisimple and any simple factor of $K_2$ has the same complexification as some simple factor of $G$. We refer to $G \times K$ as the \emph{canonical envelope} of $\Gam$.  Given $G \times K$, Margulis superrigidity theorem can be restated as:

\begin{thm}[Margulis Superrigidity Variant]
\label{thm:superrigidityvariant}
Let $G$ be a semisimple Lie group of real rank at least $2$, let $\Gamma <G$ be an irreducible lattice and $k$ a local field of characteristic zero. Then any homomorphism $\rho:\Gamma \rightarrow GL(n,k)$ extends to a homomorphism of $G \times K$ the canonical envelope of $\Gam$.
\end{thm}

We now describe a conjecture analogous to Theorem \ref{thm:superrigidityvariant} but with $\Diff(M)$ targets.
We begin by defining a \emph{local action} of a group in a way that is similar to standard definitions of pseudo-groups or groupoids but is adapted to our purposes. The definition is complicated because we need to be able to restrict local actions of topological groups to local or global actions of their discrete subgroups. So it does not suffice for our purposes to have the germ of the group action near the identity in the acting group, but rather we need it near every element in the acting group.

\begin{defn}
Let $D$ be a group and $M$ a manifold.  We say $D$ has \emph{local action} on $M$ if
\begin{enumerate}
  \item  for every point $x \in M$ and every element $d \in D$ there is an open neighorhood $V_{x,d}$ of $d \in D$ and an open neighborhood $U_{x,d}$ of $x$ in $M$ and a \emph{local action} map $\rho_{x,d} : V_{x,d} \times U_{x,d} \rightarrow M$.
  \item given $d, d'$ in $D$, whenever $z \in U_{x,d}$ and $g, hg \in V_{x,d}$ and for every $y$ such that $\rho_{x,d}(g,x) \in U_{y,d'}$ and $h \in V_{y,d'}$ we have $\rho_{x,d}(hg, z) = \rho_{y,d'}(h, \rho_{x,d}(g, z))$.
\end{enumerate}
\end{defn}

\noindent  It may be possible to offer a simpler or more transparent variant of the definition.  Point $(1)$ gives local diffeomorphisms at every point corresponding to elements of $D$, point $(2)$ requires that the collection of such local diffeomorphisms remember the group multiplication on $D$ whenever possible.  Even when $M$ is compact, one cannot restrict attention to a finite collection of local action maps unless $D$ is also compact.  The paradigmatic example to keep in mind is that $\SL(n,\R)$ acts locally on $\T^n$.  Only $\SL(n,\Z)$ has a globally defined action, but the lift to $\R^n$ one immediately has a global $\SL(n,\R)$ action which can easily be seen to give a local action on $\T^n$.
If one carries out the construction of ``blowing up" the origin in $\T^n$ as in \cite{KatokLewis}, one obtains a manifold $M$ with a local $\SL(n, \Z)$ action which does not extend to $\SL(n,\R)$ on any cover.

It is clear that one way to have a local action is to have a global action.  We say a local action \emph{restricts}
from a global action $\rho: D \times M \rightarrow M$ if we have that

$$\rho_{x,d} = \rho|_{V_{x,d} \times U_{x,d}}$$

\noindent for every $x$ in $M$ and $d \in D$.  Given a subgroup $C<D$, and a global $C$ action $\rho_C: C \times M \rightarrow$
we say that $\rho_C$ restricts from local $D$ action if

$$\rho_{x,d}|_{(C \cap V_{x,d}) \times  U_{x,d}} = \rho|_{(V_{x,d} \cap C) \times U_{x,d}}$$

\noindent for all $x \in X$ and $d \in D$.

With these definitions in hand, we can state a general superrigidity conjecture which we believe is at the heart of the phenomena observed so far in the Zimmer program.

\begin{conjec}
\label{conjec:zimmersuperrigid}
Let $G$ be a simple Lie group of real rank at least $2$ and $\Gam <G$ a lattice. Let $G \times K$ be the canonical envelope of $\Gam$ described above. Then for any compact manifold $M$ and any homomorphism
$\rho: \Gam \rightarrow \Diff(M)$ there is a local action of $G \times K$ on $M$ that restricts to the $\Gam$ action.
\end{conjec}

\begin{rem}
\begin{enumerate}
  \item For a fixed action, one expects that the local action is trivial on a finite index subgroup of $K_1$, i.e. that the local action is one of $G \times K_1' \times K_2$ where $K'_1$ is a finite quotient of $K$. This is true for $K_1$ actions by the smooth version of the Hilbert-Smith conjecture, which has been known for some time.
  \item The example of $\SL(n,\Z)$ acting on $\T^n$ shows that one needs some notion of local action to state the conjecture.  The existence of isometric actions that extend to $K_2$ justify the need for the compact extension of $G$.  The existence of actions through finite quotients of $\Gam$ justify the need for $K_1$.
\end{enumerate}
\end{rem}



At the moment, Conjecture \ref{conjec:zimmersuperrigid} incorporates all known ideas for building ``exotic" actions of lattices $\Gam$ in higher rank simple Lie groups. In addition to the parabolic induction examples discussed above, there is the blow up construction introduced by Katok and Lewis which by now has several variants \cite{KatokLewis, FisherDeformation, BenvenisteFisher, FisherWhyte, KatokRH}.  The conjecture is also of a similar flavor to
a conjecture stated in various forms by Labourie, Margulis and Zimmer that a manifold admitting a higher rank lattice action should, under some circumstances, be homogeneous on an open dense set \cite{LabICM, MargulisPCR}

In complete generality Conjecture \ref{conjec:zimmersuperrigid} seems very far out of reach.  It does seem most accessible for the case where $G=\SL(n,\R)$ and $\dim(M)=n$ and perhaps when the action is analytic.  By the work of Brown, the author and Hurtado and work of Brown, Rodriguez Hertz and Wang, the conjecture is known for lattices in $\SL(n,\R)$ for $\dim(M) < n$.  A very interesting and overlooked paper by Uchida from 1979 classifies all analytic actions of $\SL(n,\R)$ on the sphere $S^n$ \cite{Uchida}.  This suggests starting with the following

\begin{prob}
Classify analytic $\SL(n,\R)$ actions on manifolds of dimension $n$.  Classify analytic local actions of $\SL(n,\R)$ on manifolds of dimension $n$.
\end{prob}

\noindent The second part of the problem is clearly harder than the first. For both parts, it should be be useful to look at \cite{CairnsGhys, Stuck88}.

Other contexts in which Conjecture \ref{conjec:zimmersuperrigid} might be more accessible is when one assumes additional geometric or dynamical properties of the action.  Key contexts include Anosov actions \cite{BRHWAnnals, FisherSurvey2} and actions preserving rigid geometric structures \cite{GromovRigid,ZimmerNH}.  Both hyperbolicity of the dynamics and existence of geometric structures can be used to produce additional Lie groups acting on a manifold or at least on certain foliations, so these hypotheses should be helpful to find some kind of local action given a $\Gamma$ action.

\section{Other sources, other targets}

Another topic which we can only touch on briefly here is the generalization of Margulis' superrigidity theorem to other sources and targets.  The set of targets considered is quite often spaces of non-positive curvature, frequently without any assumption that the dimension is finite.  The set of sources is often broadened to more general locally compact groups. Since there is no good analogue of rank without linear structure, the most common assumption is that one has a locally compact, compactly generated group $G$ which is a product $G= G_1 \times \ldots \times G_k$ and that one has an irreducible lattice $\Gam<G$.  To give an indication, we state one particularly nice result due to Monod.  To do so we need to define
 a term.  We assume for the definition that $X$ is a geodesic metric space.

\begin{defn}
A subgroup $L < \Isom(X)$ is \emph{reduced} if there is no unbounded closed
convex subset $Y \subsetneq X$ such that $gY$ is finite (Hausdorff) distance from $Y$ for
all $g \in L$.
\end{defn}

Reduced is one possible geometric substitute for considering subgroups whose Zariski closure is
simple or semisimple.   We note that Monod proves other results in \cite{MonodJAMS} that require
weaker variants of this hypotheses, but these are more difficult to state.

\begin{thm}
\label{thm:reduced}
Let $\Gamma$ be an irreducible uniform lattice in a product $G = G_1 \times \cdots \times G_n$
of non-compact locally compact $\sigma$-compact groups with $n>1$. Let $H < \Isom(X)$ be a closed subgroup, where
$X$ is any complete $\CAT$ space not isometric to a finite-dimensional Euclidean
space. Let $\tau : \Gam \rightarrow H$ be a homomorphism with reduced unbounded image.
Then $\tau$ extends to a continuous homomorphism $\tilde{\tau}: G \rightarrow H$.
\end{thm}

We note that in the theorem, $X$ is not assumed to be locally compact.  We also note
that the theorem holds for non-uniform lattices with a mild assumption of \emph{square integrability}.
For a survey of earlier results, we point the reader to Burger's ICM address \cite{BurgerICM}.  In this context, we also mention that Gelander, Karlsson and Margulis have extended Monod's results to a broader class of non-positively curved spaces \cite{GelanderKarlssonMargulis}.  A key context for application of these kinds of results are lattices in isometry groups of products of trees \cite{BurgerMozesZimmer} and  to Kac-Moody groups  \cite{CapraceRemy}, which provide many examples of lattices in products of locally compact compactly generated groups.

A major difference between existing geometric superrigidity theorems like Monod's and Margulis' Theorem
\ref{thm:superrigidity} is that Margulis does not need any assumption like \emph{reduced}.  It is a major
open problem in the area to prove some analogue of this fact.  We state here a version of this question.
Since $X$ is not locally compact, we need to modify the notion of a representation \emph{almost extending} slightly.
If $\Gamma <G$ is a lattice and $X$ is a non-positively curved space, we say $\rho: \Gamma \rightarrow \Isom(X)$
\emph{almost extends} if there exists $\rho_1: G \rightarrow \Isom(X)$ and $\rho_2: \Gamma \rightarrow \Isom(X)$ where
$\rho_2$ has bounded image, $\rho_1(G)$ commutes with $\rho_2(\Gamma)$ and $\rho(\gamma) = \rho_1(\gamma)\rho_2(\gamma)$ for all $\gamma \in \Gam$.

\begin{qtn}
\label{qtn:geometricalmostextends}
Let $G = G_1 \times \cdots \times G_n$ for $n>2$ where each $G_i$ is a locally compact group with Kazhdan's Property $(T)$ or let $G$ be a simple Lie group of higher real rank.  Assume $\Gam \in G$ is a cocompact lattice and that
$X$ is a $\CAT(0)$ space.  Given $\rho: \Gam \rightarrow \Isom(X)$, is there a $\Gamma$-invariant subspace $Z \subset X$ such that $\rho: \Gam \rightarrow \Isom(Z)$ almost extends to $G$?
\end{qtn}

\noindent One can easily see that the passage from $X$ to $Z$ is necessary by taking a the $G$ action on $G/K$, the symmetric space for $G$, restricting to the $\Gamma$ action and adding a  discrete $\Gamma$ periodic family of rays to $G/K$.  One might assume something weaker than Property $(T)$ for each $G_i$ and should not really require the lattice is cocompact, but solving the question as formulated above would be a good first step.

Recently Bader and Furman have deeply rethought the proof of Margulis superrigidity theorem \cite{BaderFurman1,BF2,BF3}.  This work is used in the proof of Theorems \ref{thm:bfms1} and \ref{thm:charithmetic}.  It was also used quite strikingly in a proof of superrigidity theorems for groups which are not lattices in any locally compact group.  The groups in question are isometry groups of the so-called exotic $\tilde{A}_2$ buildings. The isometry groups of this buildings are known to be, in many cases, discrete and cocompact. Bader, Caprace and Lecureaux prove a superrigidity theorem for large enough groups of isometries of buildings of type $A_2$ and use this to show that a lattice in the isometry groups of a building of type $A_2$ has an infinite image linear representation if and only if the building is classical and so the isometry group is a linear group over a totally disconnected local field \cite{BCL}.

\section{The normal subgroup theorem, commensurators, attempts at unification}
\label{catchall}

In this section we describe some results and questions related to Theorems \ref{thm:commsuperrigidity} and \ref{thm:nst}.  We also describe some attempts to unify the phenomena behind Theorem \ref{thm:superrigidity} and \ref{thm:nst}.

As described in \cite{BRHWLocal} in this volume, Margulis proof of the normal subgroups theorem follows a remarkable strategy.  He proves that given a higher rank lattice $\Gamma$ and a normal subgroup $N$, the quotient group $\Gamma/N$ has Property $(T)$ and is amenable.  From this one trivially deduces that $\Gamma/N$ is a finite subgroup.  For more discussion see the article of Brown, Rodriguez Hertz and Wang in this volume.

We begin by mentioning that the normal subgroup theorem has also been generalized to contexts of products of fairly arbitrary locally compact, compactly generated groups.  This was first done by Burger and Mozes in the special
case where $G= G_1 \times \cdots \times G_k$ where each $G_i$ is a large enough subgroup of $\Aut(T_i)$ where $T_i$ is a regular tree.  Burger and Mozes used this result in order to show that certain irreducible lattices they construct in such $G$ are infinite simple groups \cite{BurgerMozes1, BurgerMozes2}.  These new simple lattices are (a) finitely presented (b) torsion-free, (c) fundamental groups of finite, locally $\CAT$-complexes, (c) of cohomological dimension 2, (d) biautomatic and (e) the free product of two isomorphic free groups $F_1$ and $F_2$ over a common finite index subgroup. The existence of such simple groups is quite surprising.  In later work, Bader and Shalom proved a much more general result about normal subgroups of lattices in fairly arbitrary products of locally compact second countable compactly generated groups \cite{BaderShalom}.  To be clear, the Bader-Shalom paper gives the ``amenability half" of the proof, i.e. that $\Gamma/N$ is amenable.  Shalom had proven earlier that $\Gamma/N$ has property $(T)$ in \cite{ShalomInv}. These results were used by Caprace and Remy to show that certain Kac-Moody groups are also simple groups \cite{CapraceRemy}.

We note that uniform lattices in rank one simple Lie groups are hyperbolic.  This means, in particular, that they
have infinitely many, infinite, infinite index normal subgroups by Gromov's geometric variants on small
cancellation theory.  In particular, the normal closure $N$ of any large enough element $\gamma \in \Gam$ has
the property that $\Gam/N$ is an infinite hyperbolic group \cite{Gromov-hyp}.  The following interesting
question is open.

\begin{qtn}
\label{qtn:rank1nst}
Let $G$ be rank $1$ simple Lie group and $\Gam <G$ a lattice.  Assume $N$ is a finitely generated infinite normal
subgroup of $\Gam$.  If $G$ is $\Sp(n,1)$ of $F_4^{-20}$, is $N$ necessarily finite index?  If $G$ is $\SU(n,1)$ or $\SU(n,1)$ is $\Gam/N$ necessarily a-$(T)$-menable?
\end{qtn}

\noindent
A-$(T)$-menability is a strong negation of Property $(T)$ introduced by Gromov.  One way to prove the question would be
to prove that for all $G$ and all $\Gamma$ and $N$, the group $\Gamma/N$ is a-$(T)$-menable.  This would  resemble Margulis' proof of the normal subgroup theorem where a key step is proving the quotient group is amenable.
For both $\SO(n,1)$ with $n>2$ and for all $n$ many lattices are known to have finitely generated normal subgroups $N$ where $\Gam/N$ is abelian, see e.g. \cite{AgolFiber, Kielak}. For $\SU(n,1)$ both abelian groups and surface groups are known to occur \cite{KapovitchSU, Stover}. For $\Gam < \SO(2,1)$ it is relatively elementary that there are no infinite index finitely generated normal subgroups.  The author first learned a variant of this question in around $2006$ from Farb.

An older question related to Theorem \ref{thm:nst} was raised in conversation between Zimmer and Margulis in the late 1970's. Given a lattice $\Gam < G$, we say a subgroup $C < \Gamma$ is \emph{commensurated} if $\Gam < \Comm_G(C)$.

\begin{conjec}
\label{qtn:marguliszimmer}
Let $G$ be a simple Lie group of real rank at least $2$ and $\Gam<G$ a lattice.  Let $N$
be a commensurated subgroup of $\Gam$.  Then $N$ is either finite or finite index in $\Gam$.
\end{conjec}

For a fairly large set of non-uniform lattices the conjecture is known by work of Venkataramana and Shalom-Willis
\cite{ShalomWillis, Venkataramana}. Shalom and Willis also formulate a natural generalization for irreducible lattices in products including $S$ arithmetic lattices. We do not include it here in the interest of brevity.

It has been known since the conversation between Margulis and Zimmer that Conjecture \ref{qtn:marguliszimmer} can be formulated as a question about homomorphisms from $\Gam$ to a certain locally compact group that is a kind of completion of $\Gam$ with respect to $N$.  Shalom and Willis prove their results on Conjecture \ref{qtn:marguliszimmer} by proving a superrigidity theorem for homomorphisms of a special class of lattices to general locally compact groups.  They also formulate an intriguing superrigidity conjecture for homomorphisms
from any higher rank lattice $\Gam$ to locally compact groups.  They demonstrate that their conjecture implies not only Conjecture \ref{qtn:marguliszimmer} but also Theorems \ref{thm:superrigidity} and \ref{thm:nst} and also the Congruence Subgroup Conjecture; see \cite[Conjecture 7.7]{ShalomWillis} and the surrounding discussion.

In the context of that work, Shalom raised a question about an interesting analogue of Corollary \ref{cor:arithmeticity}.

\begin{qtn}
\label{qtn:Shalom}
Let $G$ be a simple Lie group and $\Gam$ a Zariski dense discrete subgroup.  Assuming $\Comm_G(\Gam)$ is
not discrete, is $\Gam$ an arithmetic lattice?
\end{qtn}

Since $\Gam < \Comm_G(\Gam)$, it is relatively easy to see that simplicity of $G$ implies this is equivalent to assuming $\Comm_G(\Gam)$ is dense in $G$. For finitely generated subgroups of $\SO(3,1)$ the question was answered by Mj building on work of Leinenger, Long and Reid \cite{Mj, LLR}. For finitely generated subgroups in $\SO(2,1)$ the conjecture is easily resolved by noting the limit set is proper closed subset of $S^1$.  Mj also shows that in general it suffices to consider the case where the limit set is full.  Recent work of Koberda and Mj studies the case where there is an arithmetic lattice $\Gam_0$ such that $\Gam \lhd \Gam_0$ and resolves this case in many settings, including when $\Gam_0/\Gam$ is abelian \cite{KoberdaMj1,KoberdaMj2}.   Ongoing work of the author, Koberds, Mj and van Limbeek
seems likely to completely resolve this question \cite{FKMvL}.

Another very interesting variant of the normal subgroup theorem was raised recently by Margulis in response to the proof of results of Abert et al in \cite{SevenSamuraiLong}.

\begin{conjec}
\label{con:margulisinj}
Let $G$ be a simple Lie group of real rank at least $2$ and $\Gam <G$ a discrete subgroup.  Further
assume that the injectivity radius is bounded on $K \backslash G /\Gam$.  Then $\Gam$ is a lattice
in $G$.
\end{conjec}

It is easy to see that this conjecture implies the normal subgroup theorem.  The results in \cite{SevenSamuraiLong}
are proven using a theorem of Stuck and Zimmer, which itself is proven by using elements of Margulis' proof of the normal subgroup theorem see \cite{StuckZimmer, ZimmerIHESstabilizer}.  The conjecture above does not seem immediately accessible by variants of Margulis' proof and will probably require a new approach.

The work of Stuck and Zimmer was the precursor of a long sequence of works by Nevo and Zimmer concerning actions of higher rank simple groups.  Their most striking result is:

\begin{thm}[Nevo-Zimmer]
\label{thm:nevozimmer}
Let $G$ be a simple Lie group of real rank at least two. Let $\mu$ be a measure on $G$ whose support generates $G$ and which is absolutely continuous with respect to Haar measure.  Assume $G$ acts on compact metric space $X$ and let $\nu$ be a $\mu$-stationary measure on $X$.  Then either $\nu$ is $G$ invariant or there exists a $\nu$-measurable $G$ equivariant map $X \rightarrow G/Q$ for $Q<G$ a proper parabolic subgroup.
\end{thm}

One can view the Nevo-Zimmer theorem as a providing complete obstructions to the existence of $G$ invariant measures in terms of \emph{projective factors}  $G/Q$.  They prove a similar but slightly more technical result for actions of lattices $\Gam <G$.   One can view this as a tool for studying actions of $G$ or $\Gam$ on compact manifolds $M$.  A central element of the proof of Theorem \ref{thm:bfh} is finding enough $\Gam$-invariant measures on $M$ to control growth of derivatives of the group action.  One difficulty for using Theorem \ref{thm:nevozimmer} is that it is hard to determine in practice when a measurable projective factor exists.   Another difficulty is that to control growth of derivatives one needs to control a wider class of measures than the $\mu$-stationary ones.

In the proof of Theorem \ref{thm:bfh} on Zimmer's conjecture  we use a different method of detecting invariant measures and projective factors for $\Gam$-actions on compact manifolds $M$ that is more effective for applications.  This is developed by Brown, Rodriguez-Hertz and Wang in \cite{BRHW}.  Where Nevo and Zimmer follow Margulis and study $G$-invariant $\sigma$-algebras of measurable sets on $X$ to find the projective factor, Brown, Rodriguez Hertz and Wang study invariant measures instead.  See the paper of Brown, Rodriguez Hertz and Wang in this volume for an account of how to prove Theorem \ref{thm:nst} by their methods \cite{BRHWLocal}.

The approach of Brown, Rodriguez Hertz and Wang is particularly intriguing since they have earlier used a variant of the same method in place of Zimmer's cocycle superrigidity theorem \cite{BRHWAnnals}.  The philosphy behind the approach is generally referred to as \emph{non-resonance implies invariance}.  We close this section with a brief description of this philosophy and one implementation of it.

To employ this philosophy to actions of  a lattice $\Gamma$ one always need to pass to the induced $G$ action on $(G \times M)/\Gamma$. This allows one to use the structure of $G$, namely the root data associated to a choice of Cartan subalgebra. To explain this philosophy better, I recall some basic facts.  The Cartan subgroup $A$ of $G$ is the largest subgroup diagonalizable over $\R$, the Cartan subalgebra $\mathfrak a$ is its Lie algebra.  It has been known since the work of \'{E}lie Cartan that a finite dimensional linear representation $\rho$ of $G$ is completely determined by linear functionals on $\mathfrak a$ that arise as generalized eigenvalues of the restriction of $\rho$ to $A$.  Here we use that there is always a simultaneous eigenspace decomposition for groups of commuting symmetric matrices and that this makes the eigenvalues into linear functionals.   These linear functionals are referred to as the {\em weights} of the representation.  For the adjoint representation of $G$ on it's own Lie algebra, the weights are given the special name of {\em roots}.  Corresponding to each root $\beta$ there is a unipotent subgroup $G_{\beta} <G$ called a {\em root subgroup} and it is well known that ``large enough" collections of root subgroups generate $G$.  Two linear functionals are called {\em resonant} if one is a positive multiple of the other.  Abstractly, given a $G$-action and an $A$-invariant  object $O$, one may try to associate to $O$ a class of linear functionals $\Omega$. {\em Non-resonance implies invariance} is the observation  that, given any root $\beta$ of $G$ that is not resonant to an element of $\Omega$, the object $O$ will automatically be invariant under the unipotent root group $G^{\beta}$.  If one can find enough such non-resonant roots, the object $O$ is automatically $G$-invariant.  We will illustrate this philosphy by sketching the proof of the following theorem from \cite{BRHW}.

\begin{thm}
\label{thm:BRHWinvmeasure}
Let $G$ be a simple Lie group of real rank at least $2$, let $\Gam <G$ be a lattice. Let $Q$ be a maximal parabolic in $G$ of minimal codimension.  Assume $M$ is a compact manifold and $\rho: \Gam \rightarrow \Diff(M)$ and $\dim(M)< \dim(G/Q)$.  Then $\Gam$ preserves a measure on $M$.
\end{thm}

We begin as above by inducing the action to a $G$ action on $(G \times M)/\Gamma$ and noting that $\Gam$ invariant measures on $M$ correspond exactly to $G$ invariant measures on $(G \times M)/\Gamma$.
Taking the minimal parabolic $P<G$ and using that $P$ is amenable, one finds a $P$-invariant measure $\mu$. The goal is to prove that $\mu$ is $G$-invariant.  Once $\mu$ is $G$-invariant, disintegerating $\mu$ over the map $(G \times M)/\Gamma \rightarrow G/\Gamma$ yields a $\Gamma$-invariant measure on $M$. Since the measure $\mu$ is $P$-invariant and $A<P$, $\mu$ also clearly invariant under the Cartan subgroup $A$ and so one can try to apply the philosophy that non-resonance implies invariance by associating some linear functionals to the pair $(A, \mu)$.  The linear functionals we consider are the Lyapunov exponents for the $A$-action.

More precisely we consider the  Lyapunov exponents for the restriction of the derivative of $A$ action to the subbundle $F$ of $T((G \times M)/\Gamma)$ defined by directions tangent to the $M$ fibers in that bundle over $G/\Gamma$.  We refer to this collection of linear functionals as {\em fiberwise Lyapunov exponents}. In this context \cite[Propsition 5.3]{BRHW} shows that, given an $A$-invariant measure on $X$ that projects to Haar measure on $G/\Gamma$, if a root $\beta$ of  $G$ is not resonant with any fiberwise Lyapunov exponent then the measure is invariant by the root subgroup $G_{\beta}$.  The rest of the proof is quite simple.   The stabilizer of $\mu$ contains $P$, which implies the projection of $\mu$ to $G/\Gamma$  is Haar measure, so the proposition just described applies.  The stablizer $G_{\mu}$ of $\mu$ in $G$ is a closed subgroup containing $P$. We also know that $G_{\mu}$  contains the group generated by the $G_{\beta}$ for all roots $\beta$ not resonant with any fiberwise Lyapunov exponent. We also know that the number of distinct fiberwise Lyapunov exponents is bounded by the dimension of $M$.  Since any closed subgroup of $G$ containing $P$ is parabolic, $G_{\mu}$ is parabolic.  So either $G_{\mu}=G$ or the number of resonant roots needs to be at least the dimension of $G/Q$ for $Q$ a maximal proper parabolic. This is because given any single root $\beta$ with $G_{\beta} \nless Q$ the group generated by $G_{\beta}$ and $Q$ is $G$.  Our assumption on the dimension of $M$ immediately implies there are not enough fiberwise Lyapunov exponents to produce $\dim(G/Q)$ resonant roots, so $\mu$ is $G$-invariant.

We say a few words here on why this philosophy also works to prove superrigidity type results.  One view of the proof of superrigidity, introduced by Margulis in \cite[Chapter VII]{MargulisBook} is that one starts with an $A$ invariant section of some vector bundle over $G/\Gamma$ and then proceeds to produce a finite dimensional space of sections that is $G$ invariant.  While the proof does not rely on the non-resonance condition, it should be clear that the objects considered in that proof might be amenable to an analysis like the one above.

\section{Other criteria for a subgroup to be a lattice}
\label{sec:latticecriteria}

Many of the results and conjectures discussed so far concern criteria for when a discrete subgroup $\Gam$ in $G$ is actually a lattice or even an arithmetic lattice.  We end this paper by pointing to some more theorems and questions giving criteria for Zariski dense discrete subgroups to be lattices.  One is the recent resolution  by Benoist and Miguel of a conjecture of Margulis, building on earlier work of Oh.  Another is a question of Prasad and Spatzier that can be seen as similar to the Benoist-Miguel theorem.   Finally we mention a question of Nori which is a variant of both of these phenomena and point to some other results on Nori's question by Chatterji and Venkataramana.

To state the theorem, we recall some definitions.  Let $G$ be a simple Lie group. It is possible to state a version of the theorem for $G$ semisimple as well, but we avoid this for simplicity.  A subgroup $U$ is \emph{horospherical} if it is the stable group of an element g in G,
i.e. $U := \{u \in G | \lim_{n \rightarrow \infty} g^n u g^{-n}= e\}$.  Horospherical subgroups are always nilpotent, so a lattice $\Delta < U$ is always a discrete cocompact subgroup.

\begin{thm}[Benoist-Miguel, Oh]
\label{thm:benoistmiguel}
Let $G$ be a simple Lie group of real rank at least $2$ and  $\Gam <G$ be a discrete, Zariski dense subgroup.  Assume $\Gamma$ contains a lattice $\Delta$ in some horospherical subgroup $U$ of $G$.  Then $\Gamma$ is an arithmetic lattice in $G$.
\end{thm}

This result was conjectured by Margulis, inspired by some elements of his original proof of arithmeticity for non-uniform lattices. Uniform lattices do not intersect unipotent subgroups, so the $\Gam$ appearing in the theorem is necessarily a non-uniform lattice. For many semisimple groups $G$, the result was proved by Hee Oh in her thesis and subsequent work including joint work with Benoist \cite{Oh-Thesis, OhIJM, BenoistOh1, BenoistOh2}.  Recently Benoist and Miquel have presented a simpler proof that works in full generality \cite{BenoistMiquel}.

We present conjecture of Ralf Spatzier which has also been stated elsewhere as a question by Gopal Prasad.  The conjecture is formally somewhat similar to the theorem above, but was more inspired by work on rank rigidity in differential geometry \cite{Ballman, BurnsSpatzier}.  To make the conjeture we require a definition.  Given a countable group $\Gamma$ let $A_i$  be the subset of $\Gamma$ consisting of elements  whose centralizer contains  free abelian subgroup of rank at most $i$ as a finite index subgroup.  The rank of $\Gamma$, sometimes called the Prasad-Raghnathan rank, is the minimal number $i$ such that $\Gamma = \gamma_1 A_i \cup \dots \cup \gamma_m A_i$ for some finite set $\gamma_1, \ldots, \gamma_m \in \Gamma$.  Note that any torsion free $\Gam$ has rank at least $1$.

\begin{conjec}
\label{conjecture:ralf}
Let $G$ be a simple Lie group with real rank at least $2$.   Let $\Gam$ be a Zariski dense discrete subgroup of $G$ whose rank is the real rank of $G$.  Then $\Gam$ is a lattice in $G$.
\end{conjec}

Spatzier also asked whether it was enough for $\Gam< G$ to have rank at least $2$ to be a lattice.  
It seems to be generally believed that the statement analogous to the one in the Beniost-Miguel theorem is false.  I.e. that there is an infinite covoume group $\Gam <G$ such that some maximal diagonalizable subgroup $A < G$ such that $ A \cap \Gam$ is a lattice in $A$. One can attempt to do this by taking the group generated by some lattice $\Gam_A < A$ and some other hyperbolic element $\gamma' \in G$ where $\gamma'$ and $\Gam_A$ play ping pong in an appropriate sense.  Verifying that this works seems somewhat tricky and we do not know a complete argument.

As pointed out by Chatterji and Venkataramana \cite{ChatterjiVenkataramana}, there is a more general question of Nori related to Theorem \ref{thm:benoistmiguel}.

\begin{qtn}[Nori, 1983]
If $H$ is a real algebraic subgroup of a real semi-simple algebraic group $G$, find sufficient conditions on $H$ and $G$ such that any Zariski dense discrete subgroup $\Gam$ of $G$ which intersects $H$ in a lattice in $H$, is itself a lattice in $G$.
\end{qtn}

\noindent Chatterji and Venkataramana conjecture that the answer to this question is yes in the case where $G$ is simple non-compact and $H<G$ is a proper simple non-compact subgroup.  They also prove this conjecture in many cases when the real rank of $H$ is at least $2$. Key ingredients in their proof are borrowed from Margulis' proofs of superrigidity and arithmeticity.  A case they leave open that is of particular interest is when $H= SL(2,\R)$ and $G=\SL(n, \R)$.  Even the case $n=3$ is open and open for any choice of embedding of $SL(2,\R)$ in $SL(3,\R)$.


\bibliographystyle{AWBmath}
\bibliography{bibliography}

\end{document}